\documentclass[11pt,twoside,a4paper]{article}
\usepackage{amsfonts}
\usepackage{amssymb}
\usepackage{amsmath}

\begin{document}

\newcommand{\pa}{\partial}
\newcommand{\opa}{\overline\pa}
\newcommand{\ol}{\overline }

\numberwithin{equation}{section}

\newcommand\C{\mathbb{C}}  
\newcommand\R{\mathbb{R}}
\newcommand\Z{\mathbb{Z}}
\newcommand\N{\mathbb{N}}
\newcommand\PP{\mathbb{P}}

{\LARGE \centerline{Flexible and inflexible $CR$ submanifolds}}
\vspace{0.8cm}

\centerline{\textsc{Judith Brinkschulte\footnote{Universit\"at Leipzig, Mathematisches Institut, Augustusplatz 10, D-04109 Leipzig, Germany. 
E-mail: brinkschulte@math.uni-leipzig.de}}
 and \textsc{C. Denson Hill}\footnote{Department of Mathematics, Stony Brook University, Stony Brook NY 11794, USA. E-mail: dhill@math.stonybrook.edu\\
{\bf{Key words:}} inflexible $CR$ submanifolds, deformations of $CR$ manifolds, embeddings of $CR$ manifolds \\
{\bf{2010 Mathematics Subject Classification:}} 32V30, 32V40 }}

\vspace{0.5cm}

\begin{abstract} 
In this paper we prove new embedding results for compactly supported deformations of $CR$ submanifolds of $\C^{n+d}$:
We show that if $M$ is a $2$-pseudoconcave $CR$ submanifold of type $(n,d)$ in $\C^{n+d}$, then any compactly supported $CR$ deformation stays in the space of globally $CR$ embeddable in $\C^{n+d}$ manifolds.
This improves an earlier result, where $M$ was assumed to be a quadratic $2$-pseudoconcave $CR$ submanifold of $\C^{n+d}$. We also give examples of weakly $2$-pseudoconcave $CR$ manifolds admitting compactly supported $CR$ deformations that are not even locally $CR$ embeddable.
\end{abstract}

\vspace{0.5cm}

\section{Introduction}
 
In a previous paper \cite{BH} we introduced the concept of {\it flexible} versus {\it inflexible} $CR$ submanifolds. This is related to the $CR$ embeddability of deformations of $CR$ structures. Roughly speaking a {\it flexible} 
submanifold admits a compactly supported $CR$ deformation that "pops out" of the space of globally $CR$ embeddable manifolds. On the other hand, for an {\it inflexible} $CR$ submanifold, any compactly supported $CR$ deformation stays in the space of globally $CR$ embeddable manifolds.\\

Much work has been concentrated on $CR$ manifolds $M$ of hypersurface type
which form the boundaries of strictly pseudoconvex domains. In that
situation, $M$ is {\it inflexible}  when $\dim_{CR} M \geq 2$,
and $M$ is {\it flexible}  when $\dim_{CR} M = 1$ (even without the
assumption of strict pseudoconvexity). See example 1 in section 4.\\

Even in the situation of
$\mathrm{codim}_{CR}M = 1$ (hypersurface type) it is of interest to study what happens
for split signature of the Levi form. In that hypersurface case, 1-pseudoconcavity means
that the Levi form has at least 1 negative eigenvalue, and at least 1
positive eigenvalue; 2-pseudoconcavity means that the Levi form has
at least 2 negative and at least 2 positive eigenvalues, etc. \\

And $CR$ manifolds can have higher $CR$-codimension, in which case
{\it $q$-pseudoconcavity}
 also seems to be a fruitful concept. 
It means that for every $x\in M$ and every characteristic conormal direction $\xi$ at $x$, the scalar Levi form $\mathcal{L}_x(\xi,\cdot)$ in this conormal direction has at least $q$ positive and $q$ negative eigenvalues.
 (See
section 2 for the precise definitions.)\\

The theory of pseudoconcave $CR$ manifolds was initiated approximately 25
years ago (see \cite{HN}). Since that time
it has slowly come to light that $CR$ manifolds of higher codimension arise naturally in mathematics; i.e., such manifolds
abound, but for a long time it was ignored that they have a natural $CR$
structure. Besides typical examples of quadratic $CR$ submanifolds of $\C^{n+d}$, they also arise naturally as minimal orbits for the
holomorphic action of real Lie groups on flag manifolds. These are even
homogeneous and almost always are $q$-pseudoconcave, for some $q$.
In fact in \cite{MN} the authors follow the general method initiated by N. Tanaka of investigating manifolds endowed with partial complex structures that come from Levi-Tanaka algebras which are the canonical prolongations of pseudocomplex fundamental graded Lie algebras.
A lot of explicit such examples can be found in  \cite{MN}, \cite{HN} or \cite{HN1}.\\

When $M$ is of hypersurface type, there are some hints that the 1-pseudoconcave case (Lorentzian case) and the $q$-pseudoconcave ($q \geq 2$)
differ. For example, it is in the Lorentzian signature case where it
is possible to generalize Nirenberg's example \cite{Ni} to $\dim_{CR} > 1$, as
was done in \cite{JT}. But when $q \geq 2$, that construction does not work.
Indeed in example 6 of section 4 we present a $CR$ manifold $N$, of any
CR-codimension, which is only 1-pseudoconcave (but weakly 2-pseudoconcave)
and it is {\it flexible}. This shows that our Theorem \ref{main} below
is almost optimal. However, our $N$ is not globally $CR$ embedded into Euclidean space. Therefore it remains an open problem to find a $1$-pseudoconcave $CR$ submanifold of some Euclidean space that is flexible.\\

The main result obtained in \cite{BH} was that any $2$-pseudoconcave {\it quadratic} $CR$ submanifold of type $(n,d)$ in $\mathbb{C}^{n+d}$ is inflexible. In the present paper we obtain the same result for $CR$ submanifolds that are not necessarily assumed to be quadratic. More precisely,

\newtheorem{main}{Theorem}[section]
\begin{main}   \label{main}   \ \\
Let $M$ be a $CR$ submanifold of type $(n,d)$ in $\C^{n+d}$ that is $2$-pseudoconcave. Let $(M_a, HM_a, J_a)_{\vert a\vert < a_o}$ be a compactly supported $CR$ deformation of $(M,HM,J)$.
Then, provided $a$ is sufficiently small, given any smooth $CR$ function $f: (M,HM,J)\longrightarrow \C$,  there is a $CR$ function $f_a: (M_a,HM_a,J_a)\longrightarrow \C$ such that for any given $\ell\in\mathbb{N}$, any given compact $K$ of $M$ and arbitrary small $\varepsilon >0$, one can find a $CR$ function $f_a: (M_a,HM_a,J_a)\longrightarrow \C$ such that the $\mathcal{C}^\ell$ norm of $f-f_a$ on $K$  is less than $\varepsilon$.\\
Moreover, $f_a$ can be chosen to coincide with the given $f$ outside a compact of $M$. In particular, $(M_a, HM_a, J_a)$ is $CR$ embeddable into $\C^{n+d}$ for $a$ sufficiently close to $0$.
\end{main}

\newtheorem{corr}[main]{Corollary}
\begin{corr}   \label{corr}   \ \\
Let $M$ be a 2-pseudoconcave  $CR$ submanifold of type $(n,d)$ in $\C^{n+d}$. Then $M$ is inflexible.
\end{corr}

{\it Remark.} The same result holds, with the same proof, if $\mathbb{C}^{n+d}$ is replaced by a strictly pseudoconvex domain in $\mathbb{C}^{n+d}$. We conjecture that it also holds with $\C^{n+d}$ is replaced by an $(n+d)$-dimensional Stein manifold $X$. However, our proof relies on the results from \cite{LS}; and it is not clear if these results also hold in the more general setting of Stein manifolds.\\

In the proof of Theorem \ref{main} we use an $L^2$ vanishing result obtained in \cite{LS}, which involved heavy use of integral formulas. In \cite{BH} we were able to obtain the analogous result by employing partial Fourier transform techniques, because of the quadratic nature of $M$. However in both \cite{BH} and in the present paper we also need certain
subelliptic estimates, from \cite{FK} in codimension one, and from \cite{HN}
in higher codimension. Thus, although the results of \cite{BH} were
restricted to the quadratic case, the proofs there are more self-contained, since they do not rely on the rather complicated integral
formulas upon which \cite{LS} is based.

\section{Definitions}

An abstract $CR$ manifold of type $(n,d)$ is a triple $(M, HM, J)$, where $M$ is a smooth real manifold of dimension $2n+d$, $HM$ is a subbundle of rank $2n$ of the tangent bundle $TM$, and $J: HM \rightarrow HM$ is a smooth fiber preserving bundle isomorphism with $J^2= -\mathrm{Id}$. We also require that $J$ be formally integrable; i.e. that we have
$$\lbrack T^{0,1}M,T^{0,1}M\rbrack \subset T^{0,1}M$$
where 
$$ T^{0,1}M = \lbrace X+ iJX\mid X\in \Gamma(M,HM)\rbrace \subset \Gamma(M,\mathbb{C}TM),$$
with $\Gamma$ denoting smooth sections.

The $CR$ dimension of $M$ is $n\geq 1$ and the $CR$ codimension is $d\geq 1$.\\

$M$ admits a $CR$ embedding into some complex manifold $X$ if one can find a smooth embedding $\varphi$ of $M$ into $X$ such that the induced $CR$ structure $\varphi_\ast(T^{0,1}M)$ on $\varphi(M)$ coincides with the $CR$ structure $T^{0,1}(X)\cap\mathbb{C}T(\varphi(M))$ from the ambient complex manifold $X$.\\

Let $(M,HM,J)$ be a $CR$ manifold of type $(n,d)$ globally $CR$ embedded into some complex manifold $X$. We say that $(M,HM,J)$ admits a {\it compactly supported $CR$ deformation} if there exists a family $(M_a, HM_a, J_a)_{\vert a\vert < a_o}$ of abstract $CR$ manifolds depending smoothly on a real parameter $a$, $\vert a\vert < a_o$ and converging to $(M,HM,J)$ as $a$ tends to $0$ in the usual $\mathcal{C}^\infty$ topology; we also require that $(M_a, HM_a, J_a)= (M,HM,J)$ for every $a\not = 0$ outside some compact $K$ of $M$ not depending on $a$.\\

Note that when $(M,HM,J)$ is $CR$ embedded into some complex manifold, then one can always "punch" $M$ as to obtain compactly supported $CR$ deformations (at least locally). With the exception of $n=1$ (when the formal integrability condition is always satisfied), it can be difficult, however, to find compactly supported $CR$ deformations in the absence of local $CR$ embeddability.\\

We say that $(M,HM,J)$ is a {\it flexible} $CR$ submanifold of $X$ if it admits a compactly supported $CR$ deformation $(M_a, HM_a, J_a)_{\vert a\vert < a_o}$ such that for every sufficiently small $a \not= 0$, the $CR$ structure $(M_a, HM_a, J_a)$ is not globally $CR$ embeddable into $X$. So, for example, the Heisenberg $CR$ structure $\mathbb{H}^2$ in $\C^2$ is flexible. This follows from Nirenberg's famous local nonembeddability examples \cite{Ni}, which can be interpreted as small (local) deformations of the Heisenberg
structure on $\mathbb{H}^2$. More examples will be discussed in the last section.\\

We say that $(M,HM,J)$ is an {\it inflexible} $CR$ submanifold of $X$ if it is not flexible. That means that $(M,HM,J)$ is inflexible if and only if for every compactly supported $CR$ deformation $(M_a, HM_a, J_a)_{\vert a\vert < a_o}$ of $(M, HM, J)$, the $CR$ manifold $(M_a, HM_a, J_a)$ is globally $CR$ embeddable into $X$.\\

We denote by $H^o M=\lbrace \xi\in T^\ast M\mid < X,\xi>=0, \forall X\in H_{\pi(\xi)}M\rbrace$ the {\it characteristic conormal bundle} of $M$. Here $\pi: T M \longrightarrow M$ is the natural projection. To each $\xi\in H^o_p M\setminus \lbrace 0\rbrace$, we associate the Levi form at $p$ in the codirection $\xi:$
$$\mathcal{L}_p(\xi, X) = \xi(\lbrack J\tilde X, \tilde X\rbrack )= d\tilde\xi(X,JX) \ \mathrm{for} \ X\in H_p M$$
which is Hermitian for the complex structure of $H_p M$ defined by $J$. Here $\tilde \xi$ is a section of $H^o M$ extending $\xi$ and $\tilde X$ a section of $HM$ extending $X$. \\

Following \cite{HN} $M$ is called $q$-pseudoconcave, with $0\leq q\leq\frac{n}{2}$, if for every $p\in M$ and every characteristic conormal direction $\xi\in H^o_p M\setminus \lbrace 0\rbrace$, the Levi form $\mathcal{L}_p(\xi, \cdot)$ has at least $q$ negative and $q$ positive eigenvalues.\\

For other standard definitions related to $CR$ structures we also refer the reader to \cite{HN} or \cite{HN1}.\\

\section{Proofs}

The idea  of the proof of Theorem \ref{main} is as follows: For a given $CR$ function $f$ on $M$ we want to find a $CR$ function $f_a$ on $M_a$ which is very close to the given $f$ on $M$. Therefore we want to solve the Cauchy-Riemann equations  $\opa_{M_a}u = \opa_{M_a}f $ with $u$ having compact support and the $\mathcal{C}^k$-norms of $u$ being controlled by some $\mathcal{C}^l$-norms of $\opa_{M_a}$ (uniformly with respect to $a$). Setting $f_a = f -u_a$ then gives the desired $CR$ function on $M_a$.\\

Let $M$ be as in Theorem \ref{main}, and let $B$ be a sufficiently large Euclidean ball containing  the compact $K$ that is the support of the $CR$ deformation of $M\subset\C^{n+d}$. 
Recalling that $M$ is $2$-pseudoconcave, we have the following result from
 $\cite[\mathrm{Theorem}\ 1.0.2]{LS}$:

 \newtheorem{L2}{Proposition}[section]  
\begin{L2}  \label{L2} 
\ \\ 
Let $q = n-1$ or $q=n$, and assume $f\in L^2_{n+d,q}(M\cap B)$ satisfies $\opa_M f =0$. Then there exists $u\in L^2_{n+d,q-1}(M\cap B)$ satisfying $\opa_M u = f$. 
\end{L2}

Here we are considering (unweighted) $L^2$ spaces with respect to the induced metrics from the Euclidean metric on $\mathbb{C}^{n+d}$.
By classical Hilbert space theory (see e.g. $\cite[\mathrm{Theorem}\ 1.1.2]{H}$), one deduces from Proposition \ref{L2} the following

\newtheorem{apriori}[L2]{Proposition}
\begin{apriori}   \label{apriori}   \ \\
Let $q = n-1$ or $q=n$. Then there exists a constant $C > 0$ such that 
$$\Vert u\Vert^2 \leq C (\Vert \opa_M u\Vert^2 + \Vert \opa_M^\ast u\Vert^2)$$
for all $u\in L^2_{n+d,q}(M\cap B)\cap\mathrm{Dom}(\opa_M)\cap\mathrm{Dom}(\opa_M^\ast)$. 
\end{apriori}

Next, we use again that $M$ is 2-pseudoconcave. 2-pseudoconcavity is clearly stable under smooth, small perturbations. Therefore $M_a$ is also 2-pseudoconcave for $a$ sufficiently small, and the 2 positive  resp. 2 negative eigenvalues of the Levi form in sufficiently close characteristic conormal directions can be bounded from below resp. above independent of $a$. 
Therefore one obtains a  uniform subelliptic estimate in degrees $q\in\lbrace 0,1,n-1,n\rbrace$ (by closely looking at the proofs in \cite{FK} for $d=1$ and \cite{HN} for higher codimensions):
There exists $\varepsilon > 0$ such that
 for every compact $K$ of $M$, there exists a constant $C_K > 0$ independent of $a$ such that
\begin{equation}  \label{subelliptic}
\Vert u\Vert^2_{\varepsilon} \leq C_K (\Vert \opa_{M_a} u\Vert^2 + \Vert\opa^\ast_{M_a} u\Vert^2 + \Vert u\Vert^2)
\end{equation}
for all smooth forms $u\in \mathcal{D}_K^{p,q}(M_a)$ with support contained in $K$, $0\leq p\leq n+d$, $q \in\lbrace 0,1,n-1,n\rbrace$.\\

Combining Proposition \ref{apriori} and (\ref{subelliptic}), we can establish an $L^2$ a priori estimate in degree $(n+d,n-1)$ and $(n+d,n)$, which is uniform with respect to $a$ (in the sense that the  constant involved does not depend on $a$).

\newtheorem{uniform}[L2]{Proposition}
\begin{uniform}   \label{uniform}   \ \\
There is  $a_0> 0$ and a constant $C > 0$ such that for $q\in\lbrace n-1,n\rbrace$ we have
$$\Vert u \Vert^2\leq C (\Vert \opa_{M_a} u\Vert^2 + \Vert \opa^\ast_{M_a}u\Vert^2)$$
for all $u\in L^2_{n+d,q}(M\cap B)\cap\mathrm{Dom}(\opa_{M_a})\cap\mathrm{Dom}(\opa_{M_a}^\ast)$, $\vert a\vert < a_0$.
\end{uniform}

{\it Proof.} 
Assume by contradiction that there is a sequence $\lbrace u_{a_\nu}\rbrace\in L^2_{n+d,q}(M_{a_{\nu}}\cap B)\cap\mathrm{Dom}(\opa_{M_{a_{\nu}}})\cap \mathrm{Dom}(\opa^\ast_{M_{a_{\nu}}}) $, $a_\nu \rightarrow 0$, such that 
\begin{equation}  \label{1}
\Vert u_{a_\nu}\Vert = 1,
\end{equation}
whereas
\begin{equation}  \label{2}
\Vert \opa_{M_{a_\nu}} u_{a_\nu}\Vert^2 + \Vert \opa^\ast_{M_{a_\nu}} u_{a_\nu}\Vert^2 < a_\nu.
\end{equation}

We now want to show that $\lbrace u_{a_\nu}\rbrace$ is a Cauchy sequence. \\

Remember that $M_{a_\nu}= M$ outside $K$. We now choose a slightly larger compact $K_1$ containing $K$ in its interior, and a smooth cut-off function $\chi$ such that $\chi\equiv 1$ outside $K_1$ and $\chi\equiv 0$ in a neighborhood of $K$. Since $\opa_{M_{a_\nu}}$, $\opa^\ast_{M_{a_\nu}}$ coincide with $\opa_M$, $\opa^\ast_{M}$ outside $K$, we obtain from Proposition \ref{apriori}
$$\Vert \chi u\Vert^2 \leq  C(\Vert \opa_M (\chi u)\Vert^2 + \Vert\opa^\ast_{M}(\chi u)\Vert^2 )$$
for all $u\in L^2_{n+d,q}(M_a\cap B)\cap\mathrm{Dom}(\opa_{M_{a_{\nu}}})\cap \mathrm{Dom}(\opa^\ast_{M_{a_{\nu}}})$, which implies
\begin{equation}  \label{outsideK}
\Vert \chi u\Vert^2 \leq C^\prime (\Vert \opa_M  u\Vert^2 + \Vert\opa^\ast_{M}u\Vert^2 + \int_{K_1\setminus K} \vert u\vert^2 dV )
\end{equation}
for some constant $C^\prime > 0$.\\

On the other hand, let $\eta$ be a smooth cut-off function so that $\eta\equiv 1$ in a neighborhood of $K_1$. Then $\Vert \eta u_{a_\nu} \Vert_{\varepsilon}$ is bounded by (\ref{subelliptic}), so the generalized Rellich lemma implies that the sequence $\lbrace u_{a_\nu}\rbrace$ restricted to $K_1$ is precompact in $L^2_{n+d,q}(K_1)$. Thus it is no loss of generality to asume that the restriction of $\lbrace u_{a_\nu}\rbrace$ to $K_1$ is a Cauchy sequence. But this combined with (\ref{outsideK}) implies that $\lbrace u_{a_\nu}\rbrace$ is a Cauchy sequence in $L^2_{n+d,q}(M\cap B)$.\\

 Denote by $u_0$ the limit of this sequence. From (\ref{2}) it follows that $\opa_M u_0$ and $\opa^\ast_{M}u_0$, defined in the distribution sense, both vanish. But from (\ref{1}) it also follows that $\Vert u_0\Vert = 1$. This contradicts Proposition \ref{apriori} and therefore completes the proof of the proposition. \hfill$\square$\\

By duality, we obtain from Proposition \ref{uniform} that one can solve the $\opa_{M_a}$-equation with support in $M\cap\ol B$ in degree $(0,1)$ with a uniform constant. For this, we consider an $L^2$ variant of $\opa_{M_a}$ defined in the following way: Let $u\in L^2_{p,q}(M_a\cap B)$. We say that  $u\in\mathrm{Dom}(\opa^c_{M_a})$ and $\opa_{M_a}^c u = f$ if there exists a sequence of test forms $u_j\in\mathcal{D}^{p,q}(M_a\cap B)$ such that $u_j \rightarrow u$ in $L^2$ and $\opa_{M_a} u_j \rightarrow f$ in $L^2$. \\

\newtheorem{csupport}[L2]{Proposition}
\begin{csupport}   \label{csupport}   \ \\
There is  $a_0> 0$ and a constant $C > 0$ independent of $a$ such that
for every $f\in L^2_{0,1}(M_a)$ with $\opa_{M_a} f = 0$ and $f$ compactly supported in $M\cap B$, one can find $u\in L^2_{0,0}(M_a)$ such that $\opa_{M_a}^c u = f$  and $\Vert u\Vert \leq C \Vert f\Vert$.
\end{csupport}

{\it Proof.} Consider the operator

$$
\begin{array}{cccc}
T_f : & L^2_{n+d,n} (M_a\cap B) & \longrightarrow & \mathbb{C} \\
 & \psi & \mapsto & \int_{M_a \cap B} f \wedge \varphi,
\end{array}
$$
where $\varphi \in L^2_{n+d,n-1}(M_a\cap B)$ satisfies $\opa_{M_a}\varphi = \psi$ in the weak sense and $\Vert\varphi\Vert \leq C \Vert\psi\Vert$ (such a $\varphi$ exists by Proposition \ref{uniform}). $T_f$ is well defined. Indeed, if $\opa_{M_a}\varphi = 0$, then we may apply Proposition \ref{uniform} again and conclude that there exists $h\in L^2_{n+d,n-2}(M_a\cap B)$ satisfying $\opa_{M_a} h = \varphi$. By Stokes' theorem this implies
\begin{equation}  \label{a}
\int_{M_a \cap B} f \wedge \varphi = \int_{M_a \cap B} f \wedge \opa_{M_a} h = \int_{M_a \cap B} \opa_{M_a}(f \wedge h) = 0.
\end{equation}
Note also that $T_f$ is continuous of norm $\leq C$. Using Riesz' theorem, we conclude that there exists $u\in L^2_{0,0}(M_a)$ satisfying
$$\int_{M_a \cap B} u \wedge \opa_{M_a}\varphi = T_f(\opa_{M_a}\varphi) =  \int_{M_a \cap B} f \wedge \varphi$$
for all $\varphi\in L^2_{n+d,n} (M_a\cap B)$. Let $\vartheta_a$ be the formal adjoint of $\opa_{M_a}$ on $L^2_{\cdot,\cdot}(M_a\cap B)$. It is easy to see that $\opa^c_{M_a}$ and $\vartheta_a$ are adjoint operators on $L^2_{\cdot,\cdot}(M_a\cap B)$. (\ref{a}) implies that $(u,\vartheta_a \varphi) = (f,\varphi)$ for any $\varphi\in\mathrm{Dom}(\vartheta_a)$, which is equivalent to $\opa^c_{M_a} u = f$. \hfill$\square$\\

{\it Proof of theorem \ref{main}.}

Let $f$ be a $CR$ function on $M$. Then $\opa_{M_a}f$ has compact support and tends to zero when $a$ tends to zero. Proposition \ref{csupport} implies that we can solve the equation $\opa_{M_a} u_a = \opa_{M_a} f$ with $\Vert u_a\Vert \leq C\Vert\opa_{M_a}f\Vert$ and $u_a$ supported in $M\cap \ol B$. Hence $u_a$ is as small as we wish in $L^2_{0,0}(M_a)$, provided $a$ is small enough. It is well-known that the subelliptic estimate (\ref{subelliptic}) in degree $q=0$ implies also the following:  Suppose given a compact $K^\prime\subset M_a$  and two smooth real functions $\zeta,\ \zeta_1$ with 
$\mathrm{supp}\zeta \subset\mathrm{supp}\zeta_1\subset K^\prime$ and $\zeta_1 =1$ on $\mathrm{supp}\zeta$, then
for any integer $m\in\mathbb{N}$ there exists a constant $C_{K,m}$ such that
$$\Vert \zeta u\Vert^2_{m+\varepsilon} \leq C_{K,m} (\Vert \zeta_1\opa_{M_a}u\Vert^2_m  + \Vert \zeta_1 u\Vert^2)$$
Here $\Vert\ \Vert_m$ denotes the Sobolev norm of order $m$. But then also
  the $\mathcal{C}^\ell$-norm of $u_a$ over a given compact $K^\prime\subset M_a$ can be controlled by some $\mathcal{C}^m$-norm of $\opa_{M_a}u_a = f$, and hence made small when letting $a$ tend to zero.
Setting $f_a = f- u_a$  proves the theorem.
\hfill$\square$\\

\section{Examples of flexible $CR$ submanifolds}

The aim of this section is to provide known and new examples of flexible $CR$ submanifolds.

\begin{enumerate}

\item Rossi $\cite{R}$ constructed small real analytic deformations of the standard $CR$ structure on the 3-sphere $S^3$ in $\mathbb{C}^2$, and such that the resulting abstract $CR$ structures fail to $CR$ embed globally into $\mathbb{C}^2$. Hence $S^3$ is a flexible $CR$ submanifold of $\mathbb{C}^2$.\\
This is in contrast to higher dimensions: Any strictly pseudoconvex $CR$ manifold $M$ of $CR$ dimension $n \geq 2$ is globally $CR$ embeddable into some $\mathbb{C}^N$ by \cite{BdM}. If, in addition, $M$ is the boundary of a strictly pseudoconvex domain in $\C^{n+1}$, then $M$ is inflexible. This follows from a result by \cite{T}, since in this situation we have $H^{0,1}(M)=0$.

\item Nirenberg's famous local nonembeddability examples \cite{Ni} can be interpreted as small (local) deformations of the Heisenberg
structure on $\mathbb{H}^2\subset\mathbb{C}^2$. Since the formal integrability condition for $CR$ structures is always satisfied in dimension $3$, one can use a cut-off function to make the local deformations a compactly supported deformation of the global Heisenberg group. 

\item More generally, {\it any} $3$-dimensional $CR$ submanifold is flexible. Indeed, if $M$ has a point of strict pseudoconvexity, then one can use  the local nonembeddability result of \cite{JT} to produce a small, non-locally embeddable $CR$ deformation which is compactly supported near that point.\\
If $M$ is Levi-flat, then one first makes  arbitrary small bumps near a fixed point to get points of strict pseudoconvexity and proceeds as before.

\item $S^3 \times S^3 \in \C^4$ is an example of a flexible $CR$ submanifold of codimension $2$ (because each factor is flexible). Depending on the conormal direction, its Levi-forms have signature $++$, $--$, $+0$ or $-0$. By adding more products one can obtain flexible $CR$ submanifolds of any $CR$ codimension.

\item Let $X$ be any compact Riemann surface. Then $S^3\times X$ is flexible.

\item Let $M$ be a compact $1$-pseudoconcave $CR$ submanifold of type $(2,d)$ of some complex manifold $X$, $d$ arbitrary.  Then $N = M\times \C\PP^1$ is again $1$-pseudoconcave, and even weakly $2$-pseudoconcave (the Levi form has signature $(+-0)$ in every nonzero conormal direction).  Using ideas from \cite{Hi1} we will now show that $N$ is flexible, which indicates that Theorem \ref{main} is close to being optimal.
Indeed, by Theorem 3.2  of \cite{BH1} there exists a smooth $(0,1)$-form $\omega$ on $M$ satisfying $\opa_M\omega=0$ on $M$ such that $\omega$ is not $\opa_M$-exact on any neighborhood of any point $p$ on $M$. We will use this form $\omega$ to deform the $CR$ structure on $M\times\C\PP^1$. \\

On $\C\PP^1$ we use the two standard holomorphic charts $V_+ = \C\PP^1\setminus\lbrace \infty\rbrace$ and 
$V_- = \C\PP^1\setminus\lbrace \infty\rbrace$ given by the stereographic projection, with coordinates $z_+\in V_+\simeq\C$ and $z_-\in V_-\simeq\C$, where $z_- z_+ = 1$ on $V_-\cap V_+$. Then the usual complex structure on $\C\PP^1$ is given by $\frac{\pa}{\pa\ol z_\beta}$ on $V_\beta$ for $\beta = (+,-)$. \\

Let $U$ be an open set of $M$ such that $T^{0,1}M$ is spanned over $U$ by $\ol L_1, \ol L_2$. We define $T^{0,1}N_a$ to be 
spanned over $U\times V_\beta$ by the basis
\begin{equation}  \label{definition}
 \left\{ 
\begin{aligned}
 \ol X_0 & = &\frac{\pa}{\pa\ol z_\beta} \hspace{4cm}\\
\ol X_j & = & \ol L_j + \beta a\omega(\ol L_j) z_\beta \frac{\pa}{\pa z_\beta},\ j = 1,2
\end{aligned} \right.
\end{equation}
This gives a well defined $CR$ structure on $N$. To see that the integrability condition is valid, first note that $\lbrack \ol X_0, \ol X_j\rbrack =0$ for $j=1,2$. Moreover, by assumption on $\omega$ we have
$$ 0 = \opa_M(\ol L_1,\ol L_2) = \ol L_1(\omega(\ol L_2))- \ol L_2(\omega(\ol L_1)) - \omega(\lbrack \ol L_1,\ol L_2\rbrack),$$
thus
\begin{eqnarray*}
\lbrack \ol X_1,\ol X_2\rbrack & = & \lbrack  \ol L_1, \ol L_2\rbrack + \lbrack \ol L_1, \beta a \omega(\ol L_2) z_\beta\frac{\pa}{\pa z_\beta} \rbrack + \lbrack \beta a \omega(\ol L_1) z_\beta\frac{\pa}{\pa z_\beta}, \ol L_2\rbrack \\
& = & \lbrack  \ol L_1, \ol L_2\rbrack + \beta a \big( \ol L_1(\omega(\ol L_2)) z_\beta \frac{\pa}{\pa z_\beta} - 
\ol L_2(\omega(\ol L_1)) z_\beta \frac{\pa}{\pa z_\beta}\big) \\
& = & \lbrack  \ol L_1, \ol L_2\rbrack + \beta a \omega(\lbrack \ol L_1,\ol L_2\rbrack) z_\beta \frac{\pa}{\pa z_\beta},
\end{eqnarray*}
thus $T^{0,1}N_a$ is stable under the Lie bracket.\\

 However, for $a\not= 0$, local $CR$ embeddability of $N_a$ implies the local $\opa_M$-exactness of $\omega$. The  argument follows \cite{Hi2} or \cite{Hi}. In fact the argument shows that $N_a$ is not even locally $CR$ embeddable at any point $(t^o, z^o_\beta)$ of $N_a= M\times\C\PP^1$:
Near $t^o$, $M$ is locally $CR$ embeddable into $\C^{2+d}$ with coordinate functions $\zeta_1,\ldots, \zeta_{2+d}$. We may assume that $t= (t_1,\ldots, t_{4+d}) = (\mathrm{Re}\zeta_1,\ldots,\mathrm{Re}\zeta_{2+d}, \mathrm{Im}\zeta_1,\mathrm{Im}\zeta_2)$ are real coordinates on $M$ with $t^o=0$ in these coordinates. \\

Suppose now that we have a local $CR$ embedding of $N_a$ near $(t^o, z^o_\beta)$ by $CR$ functions $u_1(t,z_\beta), u_2(t,z_\beta),\ldots, u_{3+d}(t,z_\beta)$ with $du_1\wedge\ldots\wedge du_{3+d}\not=0 $ at $(t^o,z^o_\beta)$. Then each $u_j$ is holomorphic in $z_\beta$ in view of (\ref{definition}).

 It is then not difficult to see that  $\frac{\pa u_j}{\pa z_\beta}\not= 0$ for some $j$ at $(t^o,z^o_\beta)$. By renaming, we may assume $\frac{\pa u_{3+d}}{\pa z_\beta}\not= 0$.\\

The coordinates on $(z_1,\ldots, z_{3+d})$ on  $\mathbb{C}^{3+d}$ also define $CR$ functions on $N_a$
 and $ dz_1\wedge\ldots\wedge dz_{2+d}\wedge du_{3+d} \not= 0$ at $(t^o,u^o_\beta)$. So we arrive at a new local embedding map
$$\varphi: (t,z_\beta) \mapsto ( z_1(t,z_\beta),\ldots, z_{2+d}(t,z_\beta), u_{3+d}(t,z_\beta))$$
of some neighborhood $W$ of $(t^o,z_\beta^o)$ into $\mathbb{C}^{3+d}$. $\varphi(W)$ is a piece of a real hypersurface  in $\mathbb{C}^{3+d}$. Let $w =(w_1,\ldots, w_{3+d})$ denote the coordinates in $\mathbb{C}^{3+d}$, and consider, for points on $\varphi(W)$, the function
$$F(w) = \varphi_{\ast}\big(- \lbrack \frac{\pa u_{3+d}}{\pa z_\beta}\rbrack^{-1}\big),$$
where $\varphi_{\ast}$ is the push-forward by the diffeomorphism $\varphi$ of $W$ onto $\varphi(W)$. It follows that $F$ is a $CR$ function on $\varphi(W)$; in particular, it is holomorphic in $w_{3+d}$ by the inverse mapping theorem for holomorphic functions of one variable. On $\varphi(W)$ we may define the function
\begin{equation}  \label{function}
G(w)= \int_0^{w_{3+d}} F(w_1,\ldots, w_{2+d}, \eta) d\eta
\end{equation}
by a contour integral in the $w_{3+d}$-plane. This is well defined by the open mapping theorem from one complex variable. We now pull back to get a function $g(t,z_\beta)= \varphi^\ast G$ on $V$, which is a $CR$ function there. This can be seen by replacing $F$ in (\ref{function}) by a smooth extension $\tilde F$ of $F$ off of $\varphi(W)$ such that $\opa\tilde F_{\mid\varphi(V)} =0$ and differentiating. Next we have
\begin{equation}  
\frac{\pa g}{\pa z_\beta}  = 
 F( w_1,\ldots,w_{2+d}, u_{3+d}(t,z_\beta)) \frac{\pa u_{3+d}}{\pa z_\beta}(t,z_\beta)\hspace{4cm}\\
  =  -1,
\end{equation}
so $g(t,z_\beta) = -z_\beta + \chi(t)$, where $\chi(t)$ is a smooth "constant of integration". Now the fact that $g$ is a $CR$ function implies that $\ol X_j g =0$, hence $\ol L_{j} \chi-\beta a\omega(\ol L_j) =0$ for $j=1,2$. But for $a\not= 0$, this means that there exists a neighborhood $V^\prime$ of $t^o$ on $M$ such that $\omega$ is $\opa_M$-exact on $V^\prime$. This is a contradiction to the assumption on $\omega$. Therefore for $a\not=0$, $N_a$ is not locally $CR$ embeddable on any open neighborhood of $(t^o,z_\beta^o)$ on $N_a$. \hfill$\square$\\

\end{enumerate}

\end{document}